\documentclass[preprint,review,number,1p,10pt]{paper}
\usepackage{amssymb,graphicx,enumerate}
\def\sech{\mbox{sech}}
\def\R{\mathbb{R}}

\def\Z{\mathbb{Z}}
\def\N{\mathbb{N}}

\begin{document}

\title{Discrete conservation laws and the convergence of long time simulations of the mKdV equation}

\author{C. Gorria$^1$, M.A. Alejo$^2$ and L. Vega$^2$ \\
$^1$Applied Mathematics and Statistics Department and \\
$^2$Mathematics Department, \\
The University of the Basque Country, 48080 Bilbao, Spain}

\maketitle

\begin{abstract}
Pseudospectral collocation methods and finite difference methods have been used for approximating an important family of soliton like solutions of the mKdV equation. These solutions present a structural instability which make difficult to approximate their evolution in long time intervals with enough accuracy. The standard numerical methods do not guarantee the convergence to the proper solution of the initial value problem and often fail by approaching solutions associated to different initial conditions. In this frame the numerical schemes that preserve the discrete invariants related to some conservation laws of this equation produce better results than the methods which only take care of a high consistency order. Pseudospectral spatial discretization appear as the most robust of the numerical methods, but finite difference schemes are useful in order to analyze the rule played by the conservation of the invariants in the convergence.
\end{abstract}

\textit{keyword:}
Solitons, KdV-like equations, Finite differences, Spectral collocation methods, Conservation laws



\section{Introduction}

In this paper we study from a numerical point of view the so called breather solutions of the focusing modified Korteweg-de Vries equation:
\begin{equation}
\label{mkdv}
\left\{ \begin{array}{l} \partial_t u + \partial^3_x u + 2 \partial_x (u^3) = 0,\qquad x\in \R, \\ u(x,0) = u_0(x). \end{array} \right.
\end{equation}
This equation is a canonical non-linear dispersive equation \cite{whitham74}, and therefore it appears as a good approximation of different physical problems as the motion of the curvature of some geometric flux Ref. (\cite{goldsteinpetrich91}, \cite{nakayama92}, \cite{nakayama93}), the vortex patch, ferromagnetic vortices Ref. \cite{wexler99}, fluid mechanics Ref. (\cite{helal01}, \cite{helal02}, \cite{helal06}, \cite{helal07}), traffic models, anharmonic lattices, hyperbolic surfaces, etc.

\begin{figure}[h]
\centering
\includegraphics[width=12.cm]{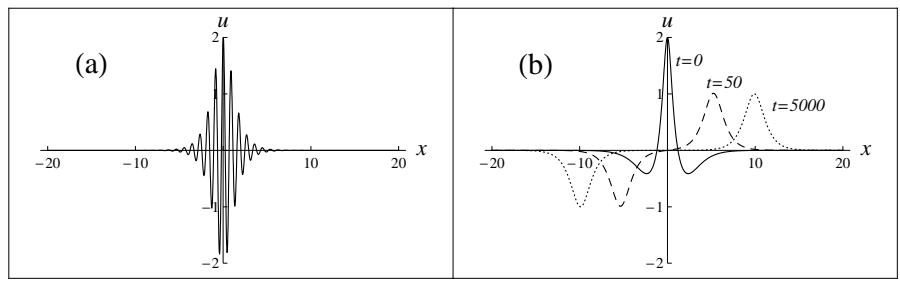}
\caption{(a) Breather type solution of mKdV at time t=0 for
$\alpha=7$ and $\beta=1$. (b) Double pole solution for $\beta = 1$
at time: $t=0$ (solid line), $t=50$ (dashed line) and $t=5000$ (doted line), after translation to the axis origin to show the logarithmic spread in time.
\label{fig1}}
\end{figure}

Travelling wave solutions that exhibit one isolated hump are easily found either by direct integration of the corresponding o.d.e. or using the inverse scattering method. In the case of the real line they are explicitly given by
\begin{equation}
\label{mkdv_soliton}
u(x,t) = \beta\sech (\beta(x-\beta^2t)),
\end{equation}
with $\beta\in\R$. Notice that the above expression could admit another two parameters to consider the traslations in space and time. As we can see, the travelling wave propagates to the right with speed $\beta^2$. From the inverse scattering point of view these solutions are characterized by the property that the reflexion coefficient has a single pole located at the imaginary axis. Solutions of more than one hump can also be constructed and they correspond to a reflexion coefficient with more than one pole in the imaginary axis. When these poles are different, they represent humps (up or down) of different heights. Therefore, these humps travel at different speed and they collide in an almost elastic way, see for example \cite{lamb80}. It is particularly interesting the degenerate case when we have only one pole, which is double. In this paper we will pay considerable attention to this situation considering solutions that evolve asymptotically in time as two equal humps such that one is up and the other is down (Figure \ref{fig1}b). It's widely used the term double pole solution to describe such solutions and we shall use it in these pages. The numerical simulations that start with this family of initial conditions converge to solutions of different type depending on the method chosen, and not all of them are satisfactory. It is of fundamental interest to clarify why some conventional numerical methods turn an initial condition from a family of solutions into a different one and which improvements in the method might prevent these instabilities.

Another relevant family of solutions is the one formed by the so called breathers (see Figure \ref{fig1}a). They were firstly obtained by M. Wadati in \cite{wadati82} and describe oscillating pulses that do not disperse. They are determined, up to translation in time and space, by two real parameters $(\alpha, \beta)$, which correspond to the frequency of the pulse and the amplitude-width of the envelope. The phase velocity of the pulse is given by $3\beta^2-\alpha^2$, and the group velocity  by $\beta^2-3\alpha^2$, $\alpha>0$, so that for $\alpha$ large with respect to $\beta$ the pulse propagates to the left with a velocity $3\alpha^2$. In fact, in this case they can be approximated by
\begin{equation}
\label{mkdv_approxbreather}
u(x,t) \approx -2\frac{\beta^2}{\alpha} \sin [\alpha(x-(3\beta^2-\alpha^2)t)]\sech[\beta(x-(\beta^2-3\alpha^2)t)].
\end{equation}

Wadati's approach to construct these solutions is based on the inverse scattering method. The breathers are characterized by the fact that the poles of the reflexion coefficient, denoted as $\pm\alpha+ i\beta$, are symmetric with respect to the imaginary axes. This family of solutions for a large $\alpha$ with respect to $\beta$ was used by Kenig, Ponce and Vega in \cite{kenig01} to prove discontinuity of the flowmaps associated to mKdV equation in the Sobolev spaces $H^s$ of functions with $s$ derivatives in $L^2(\R)$. This lack of continuity comes from the fact that two breathers with different speeds can be very close at time zero in a Sobolev norm but, because they don't disperse, they  will eventually separate and therefore the difference of their Sobolev norms becomes very big. The point is that the time of separation can be made arbitrary small by taking $\alpha$ large enough. The threshold for this lack of continuity occurs for $s=1/4$ that turns out to be sharp.

As far as we know, these solutions exist for the mKdV equation but not for the classical KdV equation (i.e. the nonlinearity $uû_x$ is changed into $uu_x$). This fact gives us a reason to study mKdV better than KdV.

A natural question that comes up from the observations we have just made, is what are the stability properties of these breather solutions when $\alpha$ is large. Due to the high oscillations of the pulse, numerical methods based on finite difference schemes do not look appropriate, and we will see later on that this is the case. However the Fourier transform of the pulse is highly concentrated around the frequency $\alpha$, and therefore pseudospectral methods seem much more natural in this case (Ref. \cite{abe80}, \cite{delahoz09} and \cite{ceniceros02}). We will see in the next sections that they are in fact very robust.

One could wonder what happens when finite difference methods are used for small $\alpha$. Notice that the solution is real analytic, and that the available well-posedness existence theory for the Cauchy problem (see for example \cite{kenig93}) allows us to conclude that the regularity is preserved along the flow and in principle, one shouldn't expect any numerical instability. Nevertheless the construction of the breather solution made by the inverse scattering method suggests that some instability can arise when we get close to the degenerate case $\alpha=0$. The reason is obvious. On one hand the double pole solution can be obtained from the two soliton solution when the poles have real parts identically zero and the imaginary parts tend to the same value. On the other hand, the same solution is obtained \cite{wadati82} by taking the limit of the breather solution when $\alpha$ goes to zero (see section 2 below for the explicit expressions). In this case the poles have an identical imaginary part while the real part is changing.

The main purpose of this paper is to see that these numerical instabilities do occur and even regular numerical methods as pseudospectral and finite differences fail when the time simulations become long enough. In the case of the finite differences, studied in section 4, two spatial discretizations have been used, with different discrete invariants associated to each. The motivation is to highlight that the numerical approximation to the two pole initial condition separates from the right solution and chooses one of the two possible branches, either two independent solitons or a breather. It is interesting to observe that the pseudospectral method, presented in section 3, is the most efficient studied here. It captures the logarithmic separation of the two humps for much longer times than the finite difference schemes. Nevertheless if a choice of a low number of harmonics or a too big step in time is made, then a wrong behavior also appears in the case of pseudospectral methods.

In the literature it is found that the accuracy of the numerical solutions is directly related with the consistency order of the numerical approach. For example the dynamical instability of the breather solutions of this equation as well as the stability of the double pole was numerically studied for reasonable long time intervals in Ref. \cite{kevrekidis04}. The convergence of the methods only guarantee good approximation in short time simulations. But these limitations may be overcome by fixing some invariant quantities, which are constant for the exact solution, along the evolution in time. These quantities are the discrete equivalent to the conservation laws of the continuous equation as the mean, the $L_2$ norm or the energy.

There are several works concerning the stability of a soliton with one hump. Also rather progress has been done in the case of multisolitons when the humps are widely separated from each other, and even more recently about the collision of two solitons, one fast and narrow and the other one slow and broad. However very little is known in the case we are considering in this paper. The available analytical techniques do not seem to apply in this case \cite{tao09}. We want to emphasize that the instability observed in this work do not lie in the arbitrary growth of the error, but in the incorrect behavior of the numerical solutions that fall within the orbits corresponding to different solutions in the phase space. In this paper we give numerical evidence that supports the difficulty of the problem.

\section{Initial condition}

We will study numerically some particular vanishing at infinity solutions of the modified Korteweg-de Vries equation (\ref{mkdv}) in 1+1 dimension, $(x,t)\in [0,T]\times \R$,
\begin{equation}
\label{boundary_cont}
\lim_{x\to \pm \infty} u(x,t) = 0.
\end{equation}

We are interested in a two parameter family of exact solutions that were first obtained by Wadati in Ref. \cite{wadati82}. They can be written as
\begin{eqnarray}
& \displaystyle u(x,t) = 2\beta \cdot\sech\left(\beta  \left(x + \gamma t\right)\right)\times \nonumber \\
& \displaystyle \frac{\cos \left[\Phi(x,t) \right] - \left(\displaystyle \frac{\beta }{\alpha}\right)\sin \left( \Phi(x,t) \right) \tanh\left(\beta(x + \gamma t)\right)}
{1+\left(\displaystyle \frac{\beta }{\alpha}\right)^2\sin^2 \left( \Phi(x,t) \right) \sech^2\left(\beta(x + \gamma t)\right)} \label{breather}
\end{eqnarray}
with $\gamma=3\alpha ^2-\beta ^2,~\delta=\alpha ^2-3 \beta ^2$ and $\Phi(x,t) = \alpha(x+\delta t)-\tan ^{-1}(\beta /\alpha)$.
They can also be written as
\begin{equation}
u(x,t) = 2\partial_x \tan^{-1} \left(\frac{G(x,t,\alpha,\beta)}{F(x,t,\alpha,\beta)}\right)
\end{equation}
with the auxiliary functions $F(x,t,\alpha,\beta)$ and $G(x,t,\alpha,\beta)$ given by
\begin{equation}
\label{f-g-aux}
\left\{ \begin{array}{l}
F(x,t,\alpha,\beta) = \cosh \left(\beta  \left(3 t \alpha ^2-t \beta ^2+x\right)\right), \\ \displaystyle
G(x,t,\alpha,\beta) = \frac{\beta  \sin \left(\alpha  \left(x+t \left(\alpha ^2-3 \beta ^2\right)\right)-\tan ^{-1}\left(\frac{\beta }{\alpha}\right)\right)}{\alpha }. \end{array} \right.
\end{equation}
Their shapes correspond to a breather with a group velocity $-\gamma=\beta^2-3\alpha^2$ when $\alpha \not= 0$ (Figure \ref{fig1}a) that degenerates to a double pole solution when $\alpha=0$ (Figure \ref{fig1}b). This second case is interesting because it behaves asymptotically as a pair of independent solitons of the mKdV Eq. (\ref{mkdv}), one facing up and the other one facing down,
\begin{equation}
\label{soliton_antisoliton}
u(x,t) = \frac{4 \beta \left(\cosh \left(\beta x-\beta^3 t\right)+\beta \left(3 \beta^2 t-x\right) \sinh
   \left(\beta x-\beta^3 t\right)\right)}{2 \left(\beta x-3 \beta^3 t\right)^2+\cosh \left(2 \beta x-2
   \beta^3 t\right)+1}.
\end{equation}
This double pole solution has three local extremes with one of them decaying to 0 asymptotically in time. The other two extremes separate with a distance $l(t)$ that becomes logarithmically large in time when $t$ is big enough, see Ref. \cite{wadati82},
\begin{equation}
\label{logarithmic_dist}
l(t) \approx \frac{2\log[4\beta^3\cdot t]}{\beta}, \qquad t >> 1.
\end{equation}
This result is easily found by taking into account the point-wise convergence of the function (\ref{soliton_antisoliton}) to 0 when time goes to $\infty$, and checking the position where the derivative of the exponential dominant terms vanishes. The property (\ref{logarithmic_dist}) will be very useful to check the accuracy of the numerical methods for $t$ big.

One of the most important properties of this type of partial differential equations is the existence of infinite conservation laws. For example the solution $u$ of the mKdV equation on the real axis $x\in \Omega=\R$ preserves the mean in space $I(u)$, the mass $\|u\|_2$ and the energy $E(u)$ over time, as detailed below. These statements are straightforward by integration by parts as well as taking into account the boundary conditions (\ref{boundary_cont}),
\begin{equation}
\label{conserv_laws}
\begin{array}{l}
\displaystyle I(u)= \int_{\Omega} u(x,t) \,\mathrm{d} x, \\
\displaystyle \|u\|_2 = \int_{\Omega} u^2(x,t) \,\mathrm{d} x, \\
\displaystyle E(u)= \int_{\Omega} \left[u^4(x,t) - u^2_x(x,t)\right] \,\mathrm{d} x. \end{array}
\end{equation}
In this paper we will show by computational experiments that in long time simulations the numerical schemes that preserve exactly some of these conservation laws, Eq. (\ref{conserv_laws}), are much more robust than the schemes that only take into account the consistency order of the method. The numerical data have been compared with the exact analytical expression for the travelling wave solutions, written by (\ref{breather}) and (\ref{soliton_antisoliton}), in order to measure its accuracy. In practice the initial condition provided by the experimental measurements, which have to be investigated numerically, may differ from those functions. The use of a powerful numerical scheme plays a crucial rule in long time experiments.

Due to the fact that the double pole initial condition (\ref{soliton_antisoliton}) is very regular and no oscillations are present, both type of methods seem to be appropriate to simulate its evolution. On the other hand, when the highly oscillating initial condition Eq. (\ref{breather}) is considered, we postulate that the pseudospectral method will be the most powerful because an acceptable number of harmonics can represent with high accuracy the shape of the solution. In this case the finite difference methods need a too high number of nodes to generate an acceptable approximation of the solution, and the global error of the numerical approximation will increase exponentially due to the big Lipschitz constant induced by the oscillations of the pulse.

\section{pseudospectral method}

As a general rule, the pseudospectral collocation methods, Refs. \cite{gottlieb77},  \cite{trefethen00} and \cite{guobenyu98}, are very suitable to approximate the travelling wave and breather type solutions of a partial differential equation. The oscillating profile of a function and its derivatives can be reproduced quite accurately by a manageable set of harmonics Ref. \cite{fornberg77}. Due to the limitations of the practical discretization in space, the real axis $\R$ has been substituted by a long finite interval $\Omega=[-L,L]$ (the numerical simulations have been made for $L=40$) and the boundary conditions (\ref{boundary_cont}) by periodic ones,
\begin{equation}
\label{periodic_cond}
u(-L) = u(L), \quad u_x(-L) = u_x(L).
\end{equation}
These conditions suggest to use an orthogonal basis made up of the complex exponential periodic functions, $\{\phi_j(x)=e^{i \omega_j x}\}$, with $w_j = 2\pi j/(2L),~j\in\Z$. The Fourier series expansion of a function $u(x,t)\in L_2$ is defined by
\begin{equation}
\label{fourier_expansion}
u(x,t) = \sum_{j=-\infty}^{\infty} \hat{u}_j(t) \phi_j(x), \quad \mathrm{where} \quad \hat{u}_j(t) = \frac{1}{2L}\int_{-L}^{L} u(x,t) \overline{\phi_j}(x) \,\mathrm{d}x.
\end{equation}

\begin{figure}[h]
\centering
\includegraphics[width=12.cm]{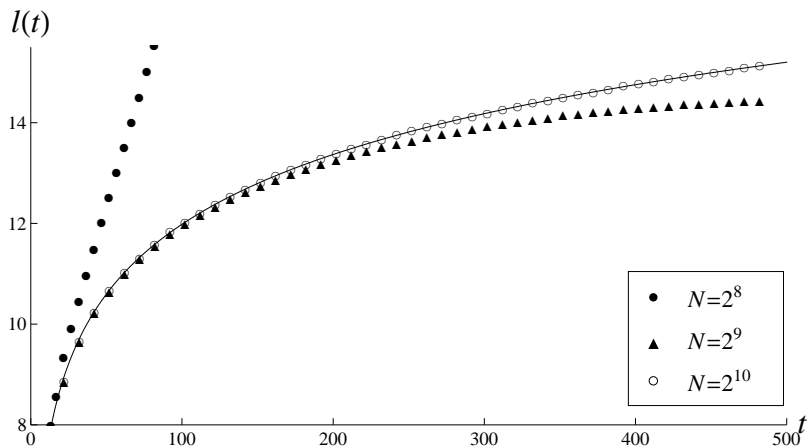}
\caption{Theoretical separation $l(t)$ of the humps of the double pole solution for $\beta = 1$ (solid line) compared with the result of pseudospectral approximation for $N=2^8$ (dots), $N=2^9$ (triangles) and $N=2^{10}$ (circles) number of collocation points and time step $\Delta t=10^{-3}$.
\label{fig4}}
\end{figure}

\begin{figure}[h]
\centering
\includegraphics[width=12.cm]{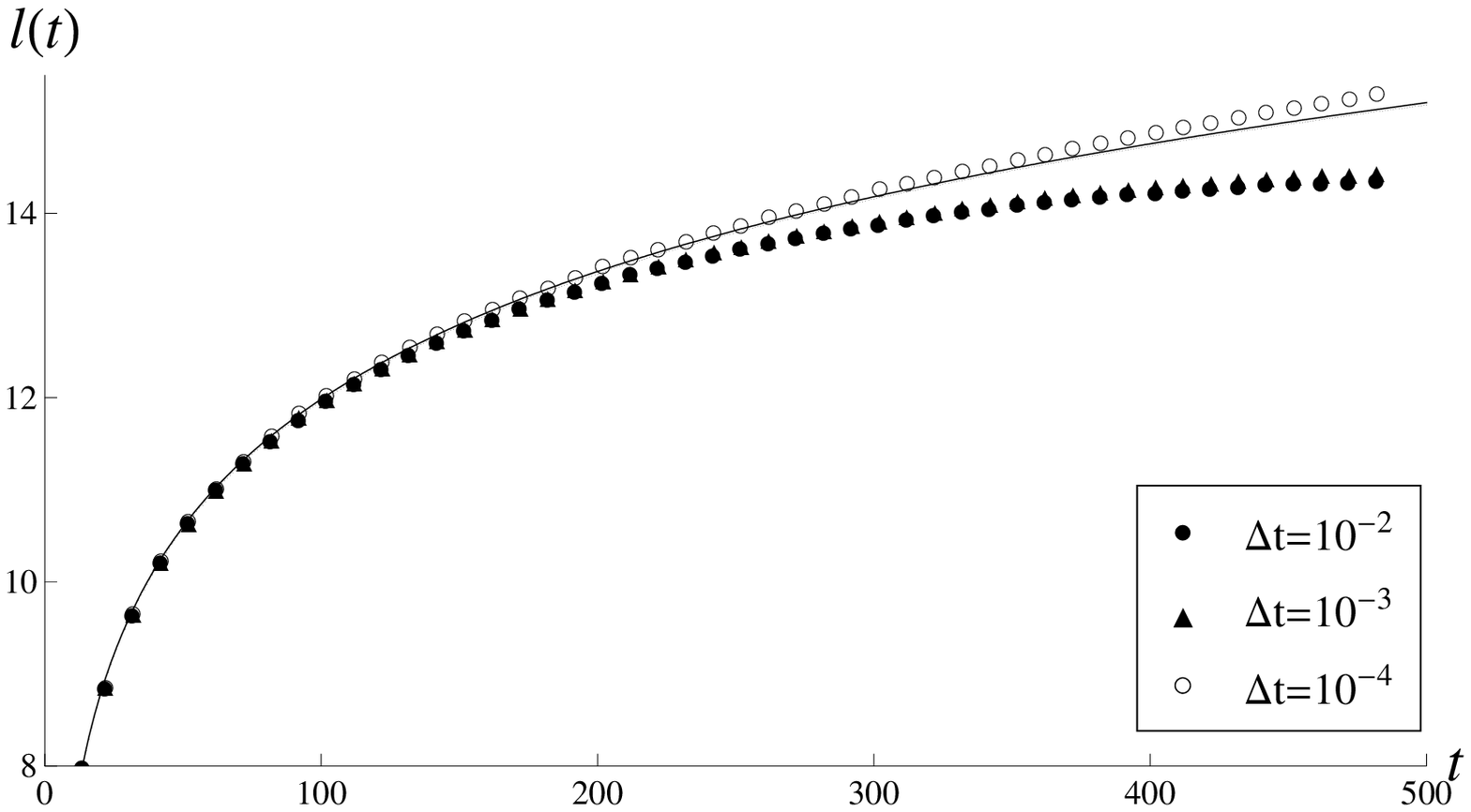}
\caption{The same as Fig \ref{fig4} for $\beta = 1$, $N=2^9$ and different time step sizes, $\Delta t=10^{-2}$ (dots), $\Delta t=10^{-3}$ (triangles) and $\Delta t=10^{-4}$ (circles).
\label{fig5}}
\end{figure}

At initial time, $t=0$, the infinite expression (\ref{fourier_expansion}) is replaced by some truncated series at the collocation points, $x_n = 2nL/N,\quad n=-N/2,\ldots,N/2-1$, to adapt it to the discrete numerical scheme,
\begin{equation}
\label{truncated_fourier}
U^{(0)}_n=u_0(x_n) = \sum_{j=-N/2}^{N/2-1} \hat{U}_j^{(0)} \phi_j(x_n), \qquad \hat{U}_j^{(0)} = \frac{1}{N}\sum_{n=-N/2}^{N/2-1} u_0(x_n) \overline{\phi_j}(x_n)
\end{equation}

The advantage of this approximation is that it can be notoriously accelerated by the well known FFT algorithm when $N=2^q,~q\in \N$, collocation points are taken. In our work we have ran simulations for $N=2^8$, $N=2^9$ and $N=2^{10}$ points to check the precision of the discretization and the robustness of the method. These calculations have been computed by a FORTRAN subroutine provided by Netlib repository \cite{gardner95}. Some authors \cite{beylkin98,du05,musluerbay03} separates the linear and the nonlinear parts of the right hand side of the equation (\ref{mkdv}) and later on they introduce an exponential integrant factor in the solution to deal with the linear part. The advantage of this reduction is that the only term to be discretized by Fourier transform or by finite differences is the nonlinear one and its norm has order $O(N)$ instead of $O(N^3)$. In this context the mentioned strategy would not produce significant improvement in the results because the wave-number of the functions of the orthogonal basis used to approximate the solutions is not necessary high. The situation is radically different when solutions that involve high frequency radiations are investigated, then it is really interesting to avoid the computation of the derivatives of the linear part by using the mentioned integrant factor technique.

Let us define $U^{(0)}=\left(U^{(0)}_{-N/2},...,U^{(0)}_{N/2-1}\right)^t$ as the vector that stores the initial data, and $U^{(k)}\approx u(t_k)$ as the successive approximations of the solution at time $t_k$. Then, the discrete Fourier transform can be considered as a linear operator $\hat{U}^{(k)} = \mathcal{F} U^{(k)}$ represented by a Vandermonde $N\times N$ matrix $\mathcal F$ with components $\mathcal F_{j,n}=\theta^{\left(j-\frac{N}{2}-1\right) \left(n-\frac{N}{2}-1\right)}$ and $\theta=e^{- i 2\pi/N}$.

\begin{figure}[h]
\centering
\includegraphics[width=12.cm]{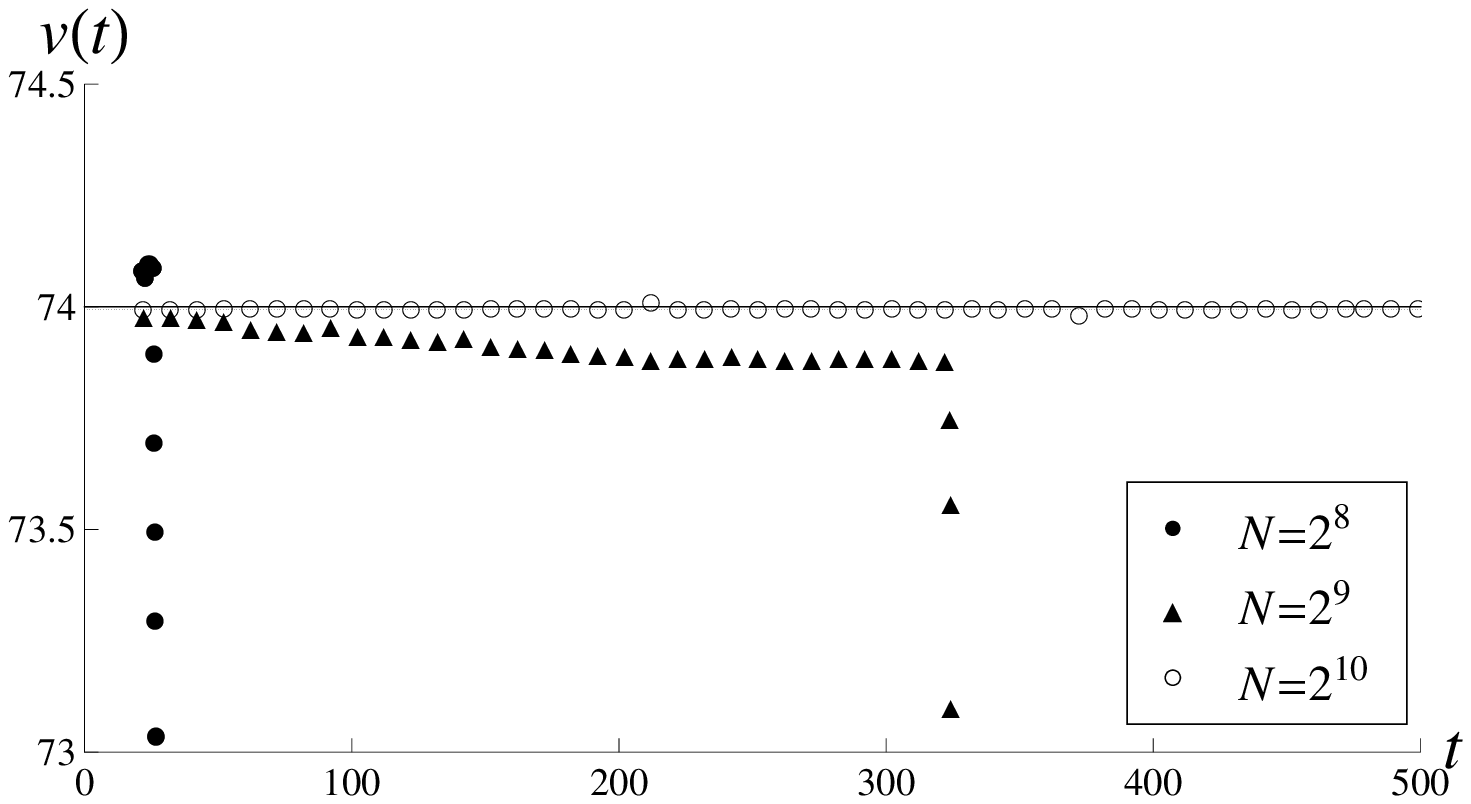}
\caption{The group velocity of the approximated breather ($\alpha = 5$ and $\beta=1$) for $\Delta t=10^{-3}$ and different number of points $N=2^8$ (dots), $N=2^9$ (triangles) and $N=2^{10}$ (circles), versus $v=3\alpha^2-\beta^2=74$.
\label{fig6}}
\end{figure}

\begin{figure}[h]
\centering
\includegraphics[width=12.cm]{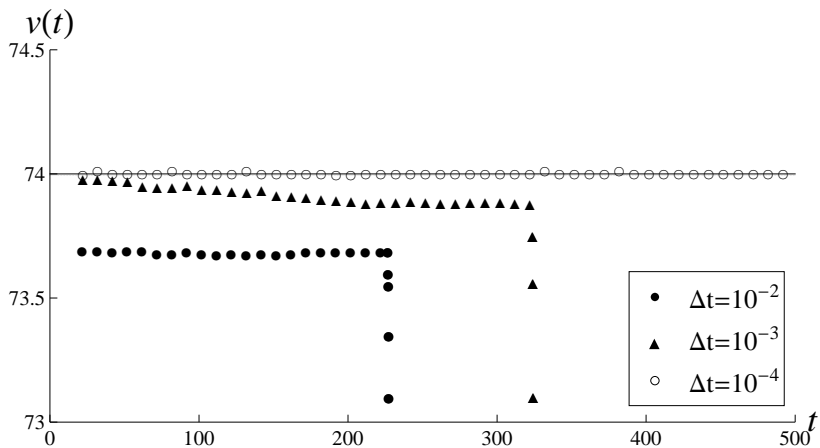}
\caption{The group velocity of the approximated breather ($\alpha = 5$ and $\beta=1$) for $N=2^9$ and different time step sizes $\Delta t=10^{-2}$ (dots), $\Delta t=10^{-3}$ (triangles) and $\Delta t=10^{-4}$ (circles), versus $v=3\alpha^2-\beta^2=74$ .
\label{fig7}}
\end{figure}

Taking into account the orthogonality $\langle\phi_j,\phi_l\rangle = \delta_{ij}$ of the periodic functions $\phi_j$ in $L_2[-\pi,\pi]$, and the expressions of their space derivatives, $\phi_j'(x) = i \omega_j \phi_j(x)$, we substitute the expansion (\ref{fourier_expansion}) in the solutions of Eq. (\ref{mkdv}) obtaining an ODE system for the Fourier coefficients,
\begin{equation}
\label{mkdv_fourier}
\left\{ \begin{array}{l} \partial_t\hat{u}_j(t) = i \left[\omega_j^3 \hat{u}_j(t) - 2 \omega_j \hat{v}_j(t)\right], \quad j=-\frac{N}{2},...,\frac{N}{2}-1, \\ \hat{u}(0) = \mathcal{F} u(x,0).  \end{array} \right.
\end{equation}
Here $\hat{v}=\mathcal{F}v$ is the Fourier transform of the nonlinear term $v=u^3$. Substituting in Eq. (\ref{mkdv_fourier}) the continuous solution $u$ and its power $v=u^3$ respectively by the $N$-dimensional vectors $U$ and $V$ with the values at the collocation points, we obtain the equivalent matrix  expression for (\ref{mkdv_fourier}),
\begin{equation}
\label{mkdv_pseudospectral}
\partial_t\mathcal{F} U = -\left[\mathcal{D}^3 \mathcal{F} U + 2\mathcal{D} \mathcal{F} V\right].
\end{equation}
The matrix $\mathcal{D}$ is diagonal with components $d_{jj}=i \omega_{j}$ and represents the spatial differentiation. Considering the linear antisymmetric operator $\mathcal{J}=\mathcal{F}^{-1}\mathcal{D}\mathcal{F}$, where  $\mathcal{J}^t=-\mathcal{J}$, the previous relation (\ref{mkdv_pseudospectral}) corresponds to the motion equation of the Hamiltonian defined by
\begin{equation}
\label{hamiltonian_pseudospectral}
\mathcal{H}_p (U) = \frac{1}{2} \left(|\mathcal{J}U|^2-\sum_{j=1}^N U_j^4 \right).
\end{equation}
Taking advantage of the conservation of the linear and the quadratic invariants Eq. (\ref{conserv_laws}) by the symplectic integrators, as far as the consistency order of the scheme, we suggest to use the implicit midpoint rule for discretization in time, together with the pseudospectral method for discretization in space. We postulate that this property about the discrete method will lead to an improvement of the time interval where the numerical solution approaches the exact one. Even though the midpoint rule is not a method with a high order of consistency, the conclusions obtained from its results will be representative and clear for the purpose of this paper. In addition the convergence for each time step is easily fulfilled, contrary to the family of Newton's type methods that are much faster than the midpoint rule but they often require deep analysis in order to avoid the problems caused by local minima.

Let $\hat{U}^{(k)}$ be the discrete Fourier transform of $U^{(k)}$, then the midpoint rule applied to the pseudospectral discretization of Eq. (\ref{mkdv_pseudospectral}) reads,
\begin{equation}
\label{pseudospectral_midpoint}
\hat{U}^{(k+1)} = \hat{U}^{(k)} - \Delta t \left(\mathcal{D}^3\hat{U}^{\left(k+\frac{1}{2}\right)} + 2\mathcal{D}\hat{V}^{\left(k+\frac{1}{2}\right)}\right),
\end{equation}
where $\hat{U}^{\left(k+\frac{1}{2}\right)}=(\hat{U}^{(k)}+\hat{U}^{(k+1)})/2$ and  $\hat{V}^{\left(k+\frac{1}{2}\right)}=(\hat{V}^{(k)}+\hat{V}^{(k+1)})/2$ are the interpolated vectors. Putting together the terms $\hat{U}^{(k+1)}$ at the left hand side of Eq. (\ref{pseudospectral_midpoint}) the norm of the subsequent functional associated to the fix-point scheme is considerably reduced. A successful implementation is provided for time steps of size $\Delta t\approx 10^{-3}$, where an acceptably cheap computational cost is involved. We get the following iteration formula
\begin{equation}
\label{pseudospectral_fixpoint}
\hat{U}^{(k+1)} = \left( \mathrm{I} + \frac{\Delta t}{2} \mathcal{D}^3\right)^{-1} \left[ \left( \mathrm{I} - \frac{\Delta t}{2} \mathcal{D}^3\right) \hat{U}^{(k)} - 2 \Delta t \mathcal{D} \hat{V}^{\left(k+\frac{1}{2}\right)} \right].
\end{equation}
Notice that the matrix $\left[ \mathrm{I} - (\Delta t/2)\mathcal{D}^3\right]$ is diagonal and its inversion is directly made. The real and the imaginary part of the system (\ref{pseudospectral_fixpoint}) are separated in different equations for the practical implementation.

In the numerical simulations with the double pole initial condition we have found a very long regular behavior of the evolution. The Figures \ref{fig4} and \ref{fig5} show the logarithmical evolution in time of the separation distance between the main humps of the approximate solutions (\ref{logarithmic_dist}) with respect to the exact $l(t)$. The increment of the number of points from $N=2^9$ to $N=2^{10}$ does not improve the accuracy as much as the reduction of the step in time from $\Delta t = 10^{-3}$ to $\Delta t = 10^{-4}$ which provides a good agreement with the exact solution. Only when the number of harmonics is too low, $N=2^8$, then the double pole initial condition can break on two independent solitons Eq. (\ref{mkdv_soliton}) as it is observed in the doted trajectory of Figure \ref{fig4}, which is almost linear and corresponds to a successive separation of the two main humps of the solution with different but nearly constant velocities.

The results for breather initial conditions with $\alpha=5$ as internal oscillation frequency have been summarized in Figures \ref{fig6} and \ref{fig7}. The choice of relaxed restrictions on the time and space steps $\Delta t > 10^{-3}$ or $N<2^9$ breaks up on the damage of the numerical solutions and the corresponding approximation to the conservation laws of the continuum equation (\ref{conserv_laws}). However a not so expensive conditions for the scheme (\ref{pseudospectral_fixpoint}) as $\Delta t \approx 10^{-4}$ or $N=2^9$ guarantee a very accurate behavior and preservation of (\ref{conserv_laws}) during long time intervals, $t \in [0,500]$, as it can be observed in Figs. \ref{fig6} and \ref{fig7}.

It is important to explain carefully the results captured on the figures (\ref{fig6}) and (\ref{fig7}) in order to avoid confusions on the interpretation. In the first one different number of collocation points have been considered while in the second one it is investigated the accuracy of the solutions for different time step sizes. In both cases the decreasing of the number of points of the mesh leads to the same effect, which is the faster loss of the convergence to the exact solution because of the jump to another family of solutions.

\section{Finite difference methods and some discrete invariants}

A wide number of papers have appeared in the last years investigating the accuracy of finite difference methods Refs. \cite{shamardan90}, \cite{helal06}, \cite{helal07} applied to nonlinear dispersive partial differential equations and comparing them with other type of discrete schemes as pseudospectral or Adomian decomposition. The goal of our analysis is to give evidences of the connection between the conservation of some invariants of the discrete scheme and the improvement on the convergence of the numerical method to the exact solution, in order to avoid the instabilities we mentioned in the introduction.

In order to implement a numerical method which approximates the solution of the initial value problem (\ref{mkdv}) we first define a spatial discrete grid of points $x_n = -L + n \Delta x,~n=1,...,N$, and the corresponding vector with the approximation to the solution in that points, $U(t)=(U_1(t),\ldots,U_N(t))$, where $U_n(t)\approx u(t,x_n)$. Let's be $\mathbf{M}(u) = -\mathbf{S}u - \mathbf{B}(u)$ the differential operator that concentrates the spatial derivatives of the mKdV equation (\ref{mkdv}), where the linear part is assumed by $\mathbf{S}(u) = \partial_x^3(u)$ and the nonlinear part by $\mathbf{B}(u) = 3u\partial_x(u^2) = 2\partial_x(u^3)$, there are several choices for an operator in finite differences $\mathbf{M_N}(\cdot) = -\mathbf{S_N}(\cdot) - \mathbf{B_N}(\cdot)$ that approximates $\mathbf{M}(\cdot)$ with different order of consistency. The time evolution of the components of the vector $U$ is governed by the following system of ODE,
\begin{equation}
\label{spatial_disc}
\partial_t U = \mathbf{M_N}(t,U),
\end{equation}
provided by the periodic boundary conditions (\ref{periodic_cond}) in $[-L,L]$, that can be written as
\begin{equation}
\label{periodic_cond_disc}
U_{n+jN} = U_{n}, \quad \forall  j\in \Z, \quad 0\le n < N.
\end{equation}

The accuracy of the numerical approximation $U$ with respect to $u$ depends on the consistency order of the operator $\mathbf{M_N}(\cdot)$, which consists of the sum of some suitable band matrices. Let's define first the forward difference matrix $\mathbf{D_1^{+}}$, where only nonzero components in each row are $d_{i,i+1}=1/\Delta x=-d_{i,i},~i=1,\ldots,N$, and $d_{N,1}=1/\Delta x$ in the corner due to periodic boundary conditions. It represents an approximation of the spatial derivative, $\partial_x u \approx \mathbf{D_1^{+}} U$, of the first order $O(\Delta x)$. Using the expression of the matrix $\mathbf{D_1^{+}}$ we can define the classical discrete approximations of the successive derivatives:

\begin{enumerate}[a)]

\item Backward differences, $\mathbf{D_1^{-}}=-(\mathbf{D_1^{+}})^T$, which approximates $\partial_x$ with order $O(\Delta x)$.

\item Central differences, $\mathbf{D_1^{c}}=(\mathbf{D_1^{+}}+\mathbf{D_1^{-}})/2$, which approximates $\partial_x$ with order $O(\Delta^2 x)$.

\item Central differences, $\mathbf{D_2^{c}}=\mathbf{D_1^{+}}\cdot \mathbf{D_1^{-}}$, which approximates $\partial^2_x$ with order $O(\Delta^2 x)$.

\item Central differences, $\mathbf{D_3^{c}}=\mathbf{D_1^{c}}\cdot \mathbf{D_2^{c}}$, which approximates $\partial^3_x$ with order $O(\Delta^2 x)$.

\end{enumerate}

It is straightforward to check that $\mathbf{D_1^{c}}$ and $\mathbf{D_3^{c}}$ are antisymmetric matrices and $\mathbf{D_2^{c}}$ is a symmetric matrix. This observation will be important to guarantee the conservation of some invariants of the discrete schemes along the time. Now we write the linear part as
\begin{equation}
\label{linear_disc}
\mathbf{S_N} U  = \mathbf{D_3^{c}}U.
\end{equation}
However the nonlinear part $\mathbf{B_N}(U)$ admits some suitable alternatives as
\begin{equation}
\label{nonlinear_disc}
\mathbf{B_N^1} (U) = 2\mathbf{D_1^{c}} U^3 \qquad \mathrm{or} \qquad \mathbf{B_N^2} (U) = 3 \mathbf{U_I} \mathbf{D_1^{c}} U^2.
\end{equation}
Here the vectors $U^2$ and $U^3$ are defined respectively by $U^2=(U_1^2,\ldots,U_N^2)$ and $U^3=(U_1^3,\ldots,U_N^3)$ and the matrix $\mathbf{U_I}$ is a diagonal matrix with the components $u_{ii}=U_i$. From now on we will suppress the $t$ dependence in the notation of $\mathbf{M_N}$ due to the autonomous structure of the mKdV equation.

Next we will choose a numerical scheme for the discretization in time, where $U^{(k)}_n \approx u(x_n,t_k)$ approximates the exact solution at time $t_k$. As we have mentioned in the previous section, the use of a symplectic integrator in the forward time step is suitable for this kind of evolution equations. They preserve the linear and quadratic invariants of the continuous equation, as far as the consistency order of the numerical method (Ref. \cite{verwer84}, [\cite{hairer02}]). Because of that we use the implicit midpoint rule that evaluates the flow of the equation $\mathbf{M_N^r}(U)= -\mathbf{S_N}U-\mathbf{B_N^r}(U),~r=1,2$, at the vector $U^{\left(k+\frac{1}{2}\right)} = (U^{(k)}+U^{(k+1)})/2$ placed between the approximations at $t_k$ and $t_{k+1}$,
\begin{equation}
\label{fd_general}
U^{(k+1)} = U^{(k)} - \Delta t \left( \mathbf{S_N}U^{\left(k+\frac{1}{2}\right)}+\mathbf{B_N^r}\left(U^{\left(k+\frac{1}{2}\right)}\right) \right), \quad r=1,2,
\end{equation}
These are finally the two implicit numerical schemes in finite differences to be analyzed.

\subsection{Convergency analysis}

The finite difference schemes proposed above will be $C$-stable (Ref. \cite{verwer84}) if the logarithmic norm of the Jacobians, $\mu_2\left(\mathbf{M_N^r}'( \overline{U})\right),~r=1,2$, are bounded independently of the grid spacing of the evolution operators for every $\overline{U}$ belonging to the segment between $U^{k+\frac{1}{2}}$ and $u(x,t+\Delta t/2)$. In both cases considered in Eq. (\ref{fd_general}) this norm is bounded by $\mu_2\left(\mathbf{M_N^r}' ( \overline{U})\right) \le \max_n (U_n^{k+\frac{1}{2}})^2$. In Ref. \cite{kenig93} it was proven that in the continuous case this is true for initial conditions $U^{0}$ in $L_{\infty}$.

To prove the $C$-stability of this implicit method we consider two different set of numerical approximations $\{U^{(k)}\}_{k\ge 0}$ and $\{\tilde{U}^{(k)}\}_{k\ge 0}$, associated to two slightly different initial conditions, and a bound for the squares of the numerical solutions $|(U_n^{\left(k+\frac{1}{2}\right)})^2|<\mu$, for all $k\ge 0$. Then
\begin{equation}
\label{c-stability1}
U^{(k+1)} - \tilde{U}^{(k+1)} = U^{(k)} - \tilde{U}^{(k)} + \Delta t \left( \mathbf{M_N^r} \left(U^{\left(k+\frac{1}{2}\right)}\right) - \mathbf{M_N^r} \left(\tilde{U}^{\left(k+\frac{1}{2}\right)}\right) \right).
\end{equation}
Now multiplying the expression (\ref{c-stability1}) by $\left(U^{\left(k+\frac{1}{2}\right)} - \tilde{U}^{\left(k+\frac{1}{2}\right)}\right)$ and applying the mean value theorem to $\mathbf{M_N^r}$ one gets that
\begin{equation}
\label{c-stability2}
\begin{array}{l}
\frac{1}{2}\|U^{(k+1)} - \tilde{U}^{(k+1)}\|^2 = \frac{1}{2}\|U^{(k)} - \tilde{U}^{(k)}\|^2 + \\
\Delta t \left( \mathbf{M_N^r} \left( U^{\left(k+\frac{1}{2}\right)}\right) - \mathbf{M_N^r} \left(\tilde{U}^{\left(k+\frac{1}{2}\right)}\right)\right) \left( U^{\left(k+\frac{1}{2}\right)} -\tilde{U}^{\left(k+\frac{1}{2}\right)} \right) \le \\
\frac{1}{2}\|U^{(k)} - \tilde{U}^{(k)}\|^2 + \Delta t \mu \|U^{\left(k+\frac{1}{2}\right)} -\tilde{U}^{\left(k+\frac{1}{2}\right)}\|^2 .\end{array}
\end{equation}
Using the triangular and the Cauchy-Schwarz inequalities we deduce
\begin{equation}
\label{c-stability3}
\|U^{(k+1)} - \tilde{U}^{(k+1)}\| \le
\|U^{(k)} - \tilde{U}^{(k)}\|^2 + \frac{\Delta t \mu}{2} \left(\|U^{(k+1)} -\tilde{U}^{(k+1)}\| + \|U^{(k)} -\tilde{U}^{(k)}\|\right)^2.
\end{equation}
To end the proof for the $C$-stability we consider the relation
\begin{equation}
\label{c-stability4}
\|U^{(k+1)} -\tilde{U}^{(k+1)}\| \, \|U^{(k)} -\tilde{U}^{(k)}\| \le \max\{\|U^{(k+1)} -\tilde{U}^{(k+1)}\|^2,\|U^{(k)} -\tilde{U}^{(k)}\|^2\},
\end{equation}
and now we can take common factor of some terms to conclude that
\begin{equation}
\label{c-stability5}
\|U^{(k+1)} - \tilde{U}^{(k+1)}\| \le
\sqrt{\frac{1 + \mu \Delta t/2}{1 - 3\mu \Delta t/2}}
\|U^{(k)} - \tilde{U}^{(k)}\|, \qquad 0 \le 3\mu \Delta t/2 < 1.
\end{equation}

The consistency of the full scheme after one step in time is proven following Ref. \cite{verwer84} defining $\hat{U}^{(k+1)} = \hat{U}^{(k+1)} + \Delta t \mathbf{M_N^r} (\tilde{U}^{\left(k+\frac{1}{2}\right)})$ as the numerical solution at time $t_{k+1}$ being $\hat{U}^{(k)}=u(t_k)$ the exact solution. Then the global error is
\begin{equation}
\label{trunc-error}
\beta(t_{k+1}) = \hat{U}^{(k+1)} - u(t_{k+1})
\end{equation}
By defining the truncation errors corresponding to the midpoint rule and the trapezoidal rule respectively as
\begin{equation}
\label{midpoint-trapezoidal-error}
\begin{array}{l}
\displaystyle d_1(t_{k+1}) = u(t_k) - u(t_{k+1}) + \Delta t \mathbf{M_N^r} \left( \frac{u(t_k) + u(t_{k+1})}{2}\right) , \\
\displaystyle d_2(t_{k+1}) = u(t_k) - u(t_{k+1}) + \frac{\Delta t}{2} \Big[ \mathbf{M_N^r} \left(u(t_k)\right) + \mathbf{M_N^r} \left(u(t_{k+1})\right)  \Big],
\end{array}
\end{equation}
then the global truncation error can be written in the following manner
\begin{equation}
\label{trunc-error2}
 \beta(t_{k+1}) = d_1(t_{k+1}) + \Delta t \left[ \mathbf{M_N^r} \left( \frac{u(t_k) + \hat{U}^{(k+1)}}{2}\right) - \mathbf{M_N^r} \left( \frac{u(t_k) + u(t_{k+1})}{2}\right) \right].
\end{equation}
Now using again the mean value theorem with the operator $\mathbf{M_N^r}$, multiplying Eq. (\ref{trunc-error2}) by $\hat{U}^{(k+1)} - u(t_{k+1})$ and using the Holder's inequality we obtain
\begin{equation}
\label{consistency}
\|\hat{U}^{(k+1)} - u(t_{k+1})\|_2 \le  \left( 1 - \frac{\mu \Delta t}{2} \right)^{-1} \|d_1(t_{k+1})\|_2, \quad \mu \Delta t <2.
\end{equation}
To complete this argument it can be deduced by Taylor expansion that, for sufficiently smooth $u$, the function $d_1$ can be bounded as $\|d_1(t_{k+1})\|_2 = O(\Delta^3 t+\Delta t\Delta^2 x)$ and consequently for $0\le \mu \Delta t < 1$
\begin{equation}
\label{consistency-bound}
\displaystyle \frac{\|\beta(t_{k+1})\|_2}{\Delta t} = O((\Delta t)^2 + (\Delta x)^2).
\end{equation}
This result finally gives us the proof of the convergence of these schemes for a given time interval $[0,T]$ whenever the steps $\Delta x$ and $\Delta t$ are chosen small enough depending on $T$. But, as we will see later, this statement does not prevent the method to corrupt the solutions for longer time intervals.

\subsection{Discrete invariants}

In the literature there are many works that emphasize the importance of designing numerical methods that preserve the invariants related to the conservation laws of the continuous model. For this purpose different strategies as geometric integration \cite{sanz-serna94} or projection methods \cite{hairer00} have been developed. In this section we investigate the transcendence of preserving some invariants by the structure of the full discretization schemes (\ref{fd_general}) for improving the stability of the numerical solution corresponding to the initial conditions, either (\ref{breather}) or (\ref{soliton_antisoliton}).

We will focus our analysis on the conservation of the following quantities, %
\begin{equation}
\label{conserv_law_disc}
\begin{array}{l}
\displaystyle \mathcal{L}_1(U) = \Delta x \sum_n U_n \approx I(u) , \\
\displaystyle \mathcal{L}_2(U) = \Delta x \sum_n (U_n)^2 = \left\langle U, U \right\rangle \approx \|u\|^2, \\
\displaystyle \mathcal{L}_3(U) = \Delta x \sum_n \left[(U_n)^4 - (U_n-U_{n-1})^2 \right] = \left\langle U^2, U^2 \right\rangle - \left\langle \mathbf{D_1^{+}} U,\mathbf{D_1^{+}} U \right\rangle  \approx E(u). \end{array}
\end{equation}
These expressions are the discrete approximations to the continuum invariants (\ref{conserv_laws}) of the mKdV equation, by the discretizations (\ref{fd_general}) presented in the previous section.
The previously proposed discrete schemes present interesting advantages with respect to the conservation of the discrete invariants (\ref{conserv_law_disc}). The first scheme, $U^{(k+1)} = U^{(k)} + \Delta t \mathbf{M_N^1}(U^{\left(k+\frac{1}{2}\right)})$, preserves the linear quantity $\mathcal{L}_1(U^k), ~\forall k\ge 0$,
\begin{equation}
\label{conserv_L1}
\mathcal{L}_1\left(U^{(k+1)}\right) = \mathcal{L}_1\left(U^{(k)}\right) + \Delta t~ \mathcal{L}_1 \left(\mathbf{M_N^1} \left( U^{\left(k+\frac{1}{2}\right)} \right)\right) = \mathcal{L}_1(U^{(k)}).
\end{equation}
Here we have used that the sum of the components of the columns of the matrixes $\mathbf{S_N}$ and $\mathbf{B_N^1}$ involved in the definition of $\mathbf{M_N^1}$ is 0, consequently $\mathcal{L}_1 (\mathbf{M_N^1} ( U^{\left(k+\frac{1}{2}\right)} )) = 0$. This property is not satisfied by the second scheme $U^{(k+1)} = U^{(k)} + \Delta t \mathbf{M_N^2}(U^{(k+\frac{1}{2})})$, however this discretization preserves $\mathcal{L}_2(U^k),$ $\forall k\ge 0$, as it can be easily shown,
\begin{equation}
\label{conserv_L2}
\begin{array}{l}
\mathcal{L}_2(U^{(k+1)}) = \left\langle U^{(k+1)}, U^{(k+1)} \right\rangle =
\left\langle U^{(k)}+ U^{(k+1)}-U^{(k)}, U^{(k)}+ U^{(k+1)}-U^{(k)} \right\rangle \\ = \mathcal{L}_2(U^{(k)}) + \left\langle U^{(k+1)} - U^{(k)}, U^{(k)} + U^{(k+1)} \right\rangle + \left\langle  U^{(k)}, U^{(k+1)} - U^{(k)} \right\rangle
\\ \displaystyle + \left\langle  U^{(k+1)} - U^{(k)}, U^{(k)} \right\rangle = \mathcal{L}_2(U^{(k)}) + \frac{\Delta t}{2} \left\langle \mathbf{M_N^2} \left(U^{\left(k+\frac{1}{2}\right)}\right), U^{\left(k+\frac{1}{2}\right)} \right\rangle = \mathcal{L}_2(U^{(k)}).
\end{array}
\end{equation}
The last simplification arises from the skew-symmetry of the matrixes $\mathbf{S_N}$ and $\mathbf{B_N^2}(U)$, which are involved in the definition of $\mathbf{M_N^2}$.

The linear and quadratic invariants shown above is all we can expect to be preserved by a symplectic method along the time steps. Their numerical precision will be of the same order as the tolerance imposed on the numerical solution of the nonlinear system (\ref{fd_general}). In addition we can prove that the exact time integration of the spatial discretization $\partial_t U = \mathbf{M_N^1}(U)$ would maintain constant along the time the invariant $\mathcal{L}_3(U)$,
\begin{equation}
\label{conserv_L3}
\begin{array}{l}
\partial_t \mathcal{L}_3(U) = \partial_t \left( \left\langle U^2, U^2 \right\rangle - \left\langle \mathbf{D_1^{+}} U,\mathbf{D_1^{+}} U \right\rangle \right) =
4\left\langle \partial_t U , U^3 \right\rangle - 2\left\langle \mathbf{D_1^{+}} \partial_t U , \mathbf{D_1^{+}} U \right\rangle \\
=
4\left\langle \partial_t U, U^3 \right\rangle + 2\left\langle (\partial_t U^T \mathbf{D_1^{-}}),\mathbf{D_1^{+}} U \right\rangle =
2\left\langle \partial_t U, \left( 2 U^3 + \mathbf{D_2^{c}} U \right) \right\rangle \\
= -2 \left\langle  \mathbf{D_3^{c}} U + 2 \mathbf{D_1^{c}} U^3 , \left( 2 U^3 + \mathbf{D_2^{c}} U \right) \right\rangle \\
= -2 \left( 2 \left\langle \mathbf{D_3^{c}} U, U^3 \right\rangle +
\left\langle \mathbf{D_3^{c}} U, \mathbf{D_2^{c}} U \right\rangle +
4 \left\langle \mathbf{D_1^{c}} U^3, U^3 \right\rangle +
2 \left\langle \mathbf{D_1^{c}} U^3, \mathbf{D_2^{c}} U \right\rangle
\right) \\
= -2 \left( 2 \left\langle \mathbf{D_3^{c}} U, U^3 \right\rangle -
\left\langle  U, \mathbf{D_3^{c}}^T \mathbf{D_2^{c}} U \right\rangle +
4 \left\langle \mathbf{D_1^{c}} U^3, U^3 \right\rangle +
- 2 \left\langle  U^3, \mathbf{D_1^{c}} \mathbf{D_2^{c}} U \right\rangle
\right) = 0
\end{array}
\end{equation}
The fact that $\mathbf{D_1^{c}}$ and $\mathbf{D_3^{c}}^T \mathbf{D_2^{c}}$ are skew-symmetric matrixes and $\mathbf{D_1^{c}}^T \mathbf{D_2^{c}} = - \mathbf{D_1^{c}} \mathbf{D_2^{c}} = -\mathbf{D_3^{c}}$ has been used to proceed with the cancelations in (\ref{conserv_L3}). The approximation to the first spatial derivative in $\mathcal{L}_3(U)$ by $\mathbf{D_1^{+}}$ and the definition of $\mathbf{M_N^1}(U)$ are crucial for this purpose. However it is easy to prove that the scheme $\partial_t U = \mathbf{M_N^2}(U)$ will not reach the condition (\ref{conserv_L3}).

Summarizing, the spatial discretization $\partial_t U = \mathbf{M_N^1}(U)$ preserves
the invariants $\mathcal{L}_1(U)$ and $\mathcal{L}_3(U)$. The first one of them is also preserved by the full discretization in time by the midpoint rule. On the other hand the full space-time discretization corresponding to $\mathbf{M_N^2}(U)$ preserves $\mathcal{L}_2(U)$. The consequences of these results on the stability of the numerical solutions of the type (\ref{soliton_antisoliton}) on long time intervals will be discussed in the next section.

\subsection{Implementation of the numerical scheme}

In this section we will describe the implementation of the schemes (\ref{fd_general}) that give approximations to the double pole solutions (\ref{soliton_antisoliton}) of the mKdV equation (\ref{mkdv}). The results will be compared with the exact solution in order to measure the efficiency of the numerical methods and to obtain conclusions.

The implicitness of the schemes makes necessary to apply them together with a numerical method for solving systems of nonlinear equations as a fix-point iteration Ref. \cite{butcher87}. First we rewrite the general method (\ref{fd_general}) as a sum of its linear and nonlinear part,
\begin{equation}
\label{fd_linear_nolinear}
U^{(k+1)} = U^{(k)} + \Delta t \mathbf{M_N^r} \left(U^{\left(k+\frac{1}{2}\right)}\right) = U^{(k)} - \Delta t \left[ \mathbf{S_N}U^{\left(k+\frac{1}{2}\right)} + \mathbf{B}_r \left(U^{\left(k+\frac{1}{2}\right)}\right) \right]  , \quad r=1,2.
\end{equation}

\begin{figure}[h]
\centering
\includegraphics[width=8.cm]{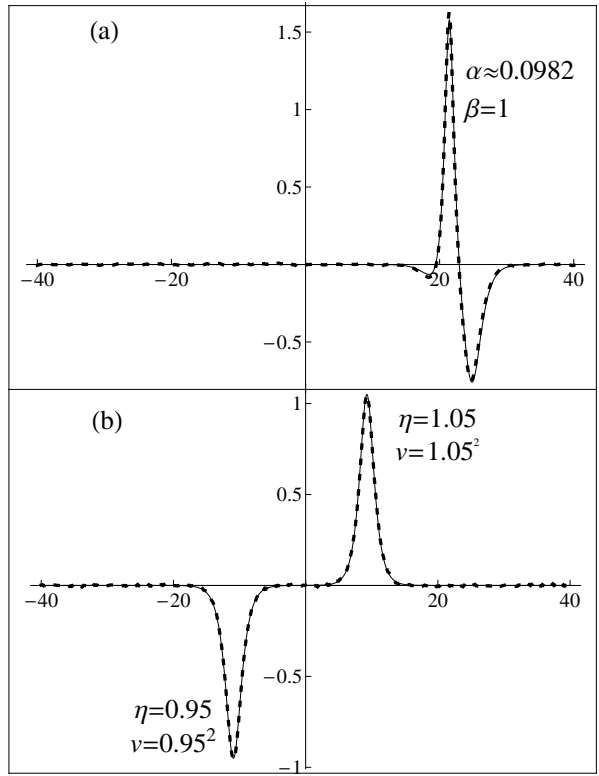}
\caption{(a) Finite difference approximation by $\mathbf{M_N^2}$ of the double pole ($\beta=1$) for $\Delta t=4\cdot10^{-2}$ and $\Delta x = 10^{-1}$ (dots) versus an exact breather of the mKdV for $\alpha=0.0982$ and $\beta=1$ (continuous line). (b) The same for $\Delta t=1\cdot10^{-2}$ versus two exact solitons of the mKdV with amplitudes $\nu=0.95$ and $\nu=1.05$.
\label{fig2}}
\end{figure}

The time step $\Delta t$ in (\ref{fd_linear_nolinear}) should be chosen small enough to assure that the operator $\Delta t \mathbf{M_N^r} (U^{k+\frac{1}{2}})$ is contractive with respect to $U^{(k+1)}$. In practice the grid size $\Delta x=1/N$ imposes a very restrictive choice of $\Delta t$ to guarantee that the norm of $\|\Delta t \mathbf{S_N}\| = O(N^3\Delta t)$ is lower than 1. This circumstance suggests to treat the implicit term $u_{xxx}$ replacing scheme (\ref{fd_linear_nolinear}) by the following formula, as was made in Ref \cite{chan85},
\begin{equation}
\label{fd_fixpoint}
U^{(k+1)} = \left( \mathbf{I}+\frac{\Delta t}{2}\mathbf{S_N}\right)^{-1} \left[ \left( \mathbf{I}-\frac{\Delta t}{2}\mathbf{S_N}\right) U^{(k)} - \Delta t \mathbf{B}_r \left(U^{(k+\frac{1}{2})}\right) \right],
\end{equation}
Now the norm of the fix-point operator is independent on $N$.

When a double pole initial condition (\ref{soliton_antisoliton}) is chosen, the results of the simulations fall into a very reach casuistic, and gives a strong hint about the instable nature of this type of solutions in the continuum case. As it was mentioned in the introduction, this particular solution is related to the fact that the reflection coefficient that appears in the inverse scattering formulation has a second order pole in the imaginary axis, Ref. \cite{wadati82}. This situation can be understood as a limiting case of solutions either made as superposition of two independent solitons (the poles are different and are locate on the imaginary axis, $\beta_1 i \approx \beta_2 i$) or as a breather of low frequency (the poles are different and symmetric with respect to the imaginary axis, $\pm \alpha + \beta i$, $\alpha < < 1$).

\begin{figure}[h]
\centering
\includegraphics[width=13.cm]{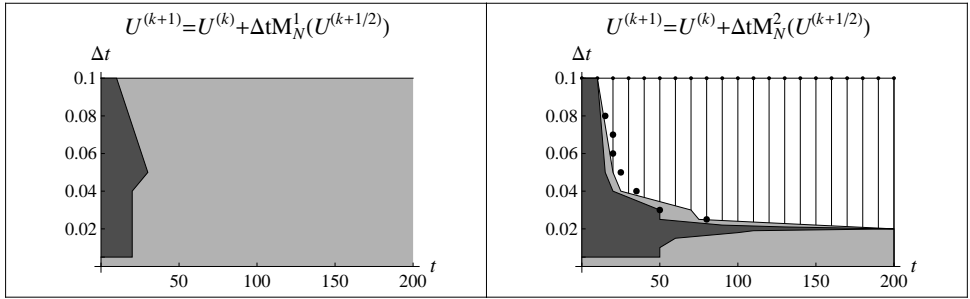}
\caption{Evolution of the solution by finite differences ($\Delta x = 10^{-1}$) for different $\Delta t$ and $\mathbf{M_N^1}$ discretization (left) or $\mathbf{M_N^2}$ (right), from a double pole (dark regions) to two independent solitons (gray regions) or a low oscillation breather (dashed region). The big dots indicates when $\mathcal{L}_3$ abruptly oscillates.
\label{fig3}}
\end{figure}

This complex scenario, where structurally different type of solutions are distinguished by the slight change on the choice of the parameters $\alpha$ and $\beta$, is captured by the behavior of the approximations produced by finite difference methods. We observed the following phenomena:

a) For long time experiments, the double pole solution is better approximated by the scheme where the spatial discretization is approached by the operator $\mathbf{M_N^2}$ than the other choice (dark region in Figure \ref{fig3}). A time step as small as $\Delta t \approx 0.02$ has to be taken to guarantee stability.

b) After reaching a transient time, the initial double pole shape is trapped by the orbits of two independent solitons when $\mathbf{M_N^1}$ spatial discretization is used or $\mathbf{M_N^2}$ with too small time step (gray region in Figure \ref{fig3}). An example of this phase change is shown in Figure \ref{fig2}(b), where $\mathbf{M_N^2}$ spatial discretization with $\Delta t = 10^{-2}$, causes that the initial double pole shape is captured by two independent solitons with amplitudes $\beta = 0.95$ and $\beta = 1.05$. A similar behavior would have been obtained for the discretization $\mathbf{M_N^1}$ with any time step after a prudential time interval $T> 50$.

c) When $\mathbf{M_N^2}$ spatial discretization is used together with not so small time step ($\Delta t > 0.02$), then the original double pole solution, after going through the shape of two independent solitons for a short period of time, leads to a breather solution with a small $\alpha$ parameter (dashed region in Figure \ref{fig3}). This jump on the shape of the solution is detected at the same time by a breakdown of the energy $\mathcal{L}_3(U)$ as is shown by the position of the big dots in the right diagram of Figure \ref{fig3}. The use of $\mathbf{M_N^1}$ discretization avoids this pathological behavior due to the intrinsic preservation of the energy. An example of the trapping by the orbit of the solution (\ref{breather}) with parameter $\alpha = 0.0982$ is shown in Figure \ref{fig2}(a).

The conclusion of this section suggests that the methods that conserve the mean $\mathcal{L}_1$ and the energy $\mathcal{L}_3$ preserve the numerical double pole solution from jumping to the orbit of a breather in the phase space of solutions. This pathological behavior is specially caused by the use of $\mathbf{M_N^2}$ discretization together with a roughly choice of time step $\Delta t$. On the other hand the discretization for $\mathbf{M_N^2}$ guarantees good results for a longer interval of time than $\mathbf{M_N^1}$ when a small enough $\Delta t$ is chosen. It means that the conservation of the discrete $\mathcal{L}_2$  is an important feature to be taken into account by a numerical method that reproduces the exact behavior of the solution for long time simulations.

\section{Conclusion}

In this paper we have studied the response of the pseudospectral methods and the finite difference methods when they are applied to a particular family of solutions of the mKdV equation. These solutions, known as breathers, depend on two parameters, $\alpha$ and $\beta$, which control respectively the internal oscillations as well as the effective group velocity of the waves. When the $\alpha$ parameter vanishes, then a special type of solution called ``double pole'' appears. Its shape is characterized by two main humps, one up and the other one down, which progressively separates with a logarithmic velocity. This solution can be considered as a limiting case of the breather solution when $\alpha \to 0$ and also it can be seen as the superposition of a soliton and an antisoliton. Therefore from the the point of view  of the inverse scattering theory the double pole solution has to be seen as a degenerate case where three type of solutions are very close. This circumstance explains that when perturbations are considered the corresponding solutions are captured by any of the three possible orbits. The first two branches are either to follow a breather with $\alpha>0$ or a soliton type solution ($\alpha=0$). In this second possibility another two orbits are possible, either the wrong branch where the amplitudes of the two humps become different and they separate by a linear speed, or the right one where the amplitudes tend to the same value and the separation is at a logarithmic rate.

The type of smooth solutions of the mKdV equation considered here clearly suggests that pseudospectral methods are the most appropriate for describing their dynamics. In fact, the functions of the basis of the Fourier discrete expansion adapt themselves very easily to the oscillations of the solution, and it is not necessary a high amount of harmonic functions to obtain an accurate approximation. However, this intrinsic oscillations of the breathers  involve a big Lipschitz constant even for $\alpha>0$ not necessarily large, and they make  the finite difference schemes not feasible to obtain good approximations in this case. Therefore we have only used them for the simulations that have the double pole as initial condition.

It is well known that the mKdV equation has an infinite number of conservation laws. As a consequence of it, we propose numerical schemes that preserve some of them with the aim of getting accurate approximations to the solutions in long time intervals. In order to increase the possibilities of success, the time advancing is reached by a symplectic method as the midpoint rule.  This implicit method involves a system of nonlinear equations that is conveniently manipulated to be economically solved by using a fix-point procedure.

The results of the simulations made by finite difference methods bring to light the relation between each family of solutions and the preservation of some of the conservation laws. Keeping constant the discrete invariants equivalent to $\mathcal{L}_1(U)$ and specially to the energy $\mathcal{L}_3(U)$. The discretization that keeps constant the equivalent to the integral of the function prevents the initial condition to degenerate from a double pole to a breather. Unfortunately after a short time the solution drops in two independent soliton-antisoliton shape. However the discretization that keeps constant the equivalent to the integral of the square of the function, $\mathcal{L}_2(U)$, reproduces for a rather longer time the double pole solution with a good accuracy. Nevertheless the approximation deteriorates for later times and jumps either to a breather type solution when the time step is rough $\Delta t>2\cdot10^{-2}$, or to a soliton-antisoliton solution when the step is small enough.

The pseudospectral method has been successfully developed with both types of initial conditions, the double pole as well as the breather. Taking $\Delta t = 10^{-4}$ and $N=2^9$ for a double pole, and $\Delta t = 10^{-5}$ and $N=2^9$ for a breather, the method reproduces with extremely good efficiency the exact solution during long intervals of time, $t\in [0,500]$. An important remark is that in these conditions, the pseudospectral discretization keeps constant the three discrete quantities investigated in this work, as far as the precision order of the method allows. Considering implementations of comparable computational time for the pseudospectral and finite difference methods, the first ones return good results for time intervals of one order of magnitude longer than the second ones.

\section{Acknowledgment}

The authors acknowledges the very helpful discussion with Dr. Ander Murua about the efficiency and the implementation of the numerical methods. The financial support of this project is acknowledged to the Spanish Ministerio de Educaci\'on y Ciencia (grant MTM2007-62186) and the Departamento de Educaci\'on, Universidades e Investigaci\'on of the Basque Government (grant IT-305-07 for research groups).

\bibliographystyle{elsarticle-num}

\end{document}